\documentclass{elsart}
\usepackage{amssymb,amsmath,amsfonts} 

\def\proofr{{\bf P\,r\,o\,o\,f\,.\,~}}

\def\proofend{$\blacktriangle$\vspace{0.4em}}

\newcommand\btheorem[1][\hspace{-1ex}]{\par\noindent\refstepcounter{thrm} \bf Theorem~\thethrm\hspace{1ex}#1.~\it}
\newcommand\etheorem{\rm\par}
\newcommand\blemma[1][\hspace{-1ex}]{\par\noindent\refstepcounter{lmm} \bf Lemma~\thelmm\hspace{1ex}#1.~\it}
\newcommand\elemma{\rm\par}
\newcommand\bpro[1][\hspace{-1ex}]{\par\noindent\refstepcounter{prp} \bf Proposition~\theprp\hspace{1ex}#1.~\it}
\newcommand\epro{\rm\par}
\newcommand\bcor[1][\hspace{-1ex}]{\par\noindent\refstepcounter{crll} \bf Corollary\hspace{1ex}#1.~\it}
\newcommand\ecor{\rm\par}
\newcommand\bnote[1][\hspace{-1ex}]{\par\noindent\refstepcounter{rmrk} \bf Remark~\thermrk\hspace{1ex}#1.~\rm}
\newcommand\enote{\rm\par}
\newcommand\bdefin[1][\hspace{-1ex}]{\par\noindent\refstepcounter{dfnt} \bf Definition~\thedfnt\hspace{1ex}#1.~\rm}
\newcommand\edefin{\rm\par}

\newcommand{\bexamp}{\begin{examp}\point\rm}
\newcommand{\Bexamp}{\begin{examp}\hspace{-1.75mm}\rm}
\newcommand{\eexamp}{\end{examp}}

\def\A{\Sigma}

\def\fromOneTo#1{[#1]}

\def\eqdf{\stackrel{\rm def}\Longleftrightarrow}
\def\defeq{\stackrel{\rm \scriptscriptstyle def}=}
\def\equivvv{\,\Leftrightarrow\,}

\journal{Discrete Mathematics}

\begin{document}

\begin{frontmatter}
\title{On reducibility of $n$-ary quasigroups}
\author{Denis S. Krotov
}
\ead{krotov@math.nsc.ru}
\address{Sobolev Institute of Mathematics, pr-t Ak. Koptyuga, 4,
Novosibirsk, 630090, Russia}
\begin{abstract}
An $n$-ary operation $Q:\Sigma^n\to \Sigma$ is called
an {$n$-ary quasigroup} {of order $|\Sigma|$} if
in the equation $x_{0}=Q(x_1, \ldots , x_n)$ knowledge of any $n$ elements
of $x_0$, \ldots , $x_n$ uniquely specifies the remaining one.
$Q$ is permutably reducible if
$Q(x_1,\ldots,x_n)=P\left(R(x_{\sigma (1)},\ldots,x_{\sigma (k)}),
\linebreak[1]
x_{\sigma (k+1)},\ldots,x_{\sigma (n)}\right)$
where $P$ and $R$ are $(n-k+1)$-ary and $k$-ary quasigroups, $\sigma$ is a permutation,
and $1<k<n$.
An $m$-ary quasigroup $S$ is called a retract of $Q$ if it
can be obtained from $Q$ or one of its inverses by fixing $n-m>0$ arguments.
We prove that if the maximum arity of a permutably irreducible retract 
of an $n$-ary qua\-si\-group $Q$ belongs to
$\{3,\ldots,n-3\}$, then $Q$ is permutably reducible.
\end{abstract}
\end{frontmatter}


\section{Introduction.}
We continue the investigation of $n$-qua\-si\-groups of order $4$ that was started in
\cite{KroPot:Lyap,Kro:2codes,PotKro:asymp}.
The general line of inquiry is the characterization of irreducible $n$-qua\-si\-groups
(which cannot be represented as a repetition-free superposition
of multary qua\-si\-groups of smaller orders).
For these reasons, we derive a new test for reducibility.
In particular, every irreducible $n$-qua\-si\-group
does not satisfy the hypothesis of the test;
this gives a new necessary condition for an $n$-qua\-si\-group
to be irreducible.
Although, historically, this work is a part of an investigation
of $n$-qua\-si\-groups of order $4$, the test,
which is given in terms of decomposability of retracts,
is suitable for any, even infinite, order.

In general, it is very natural to consider possible
representations of an $n$-qua\-si\-group
as repetition-free superpositions.
An extremely useful fact is that there exists
a unique (in some sense) canonical decomposition \cite{Cher}
(it is remarkable that this is true for essentially more wide class of functions
than the $n$-qua\-si\-groups, see \cite{Sokh:1996}).
Using the canonical decomposition of an $n$-qua\-si\-group,
it is possible to derive decompositions
for some of its retracts. The approach of this paper is opposite:
using decompositions of some retracts, we reconstruct a decomposition of
the original $n$-qua\-si\-group.

Let $\A$ be a nonempty set and $\A^n$ be the set of words of length $n$ over the alphabet $\A$.
We assume that $\A$ contains $0$; denote $\bar 0 \defeq (0,\ldots,0)$.
Let $\fromOneTo{n} \defeq \{1,\ldots,n\}$.

\bdefin[($n$-quasigroup)] An $n$-ary operation $q: \A^n \to \A$ such
that in the equality $q(x_1, \ldots , x_n) = x_{n+1}$ knowledge of any $n$ elements
of $x_1$, \ldots , $x_n$, $x_{n+1}$ uniquely specifies the remaining one is called
an \emph{$n$-ary quasigroup} of order $|\A|$ \cite{Belousov} or simply \emph{$n$-quasigroup};
we will also use the term \emph{multary quasigroup}
when the arity is not specified or inessential.
\edefin

We see that the definition is symmetric with respect to all variables
$x_1$, \ldots , $x_n$, $x_{n+1}$, while the form $q(x_1, \ldots , x_n) = x_{n+1}$ is not;
this is not handy sometimes. For this reason, we will also use the $(n+1)$-ary
predicate $q\langle \cdot \rangle$ instead:
\begin{equation}\label{e[]}
q\langle x_1, \ldots , x_n, x_{n+1}\rangle  \eqdf q(x_1, \ldots , x_n) = x_{n+1}.
\end{equation}
(In fact, the predicate $q\langle \cdot \rangle$ represents the graph of $q$.)
We use upper-case letters to name multary quasigroups in predicative form,
see the following definition for example.
It is also sometimes convenient  to talk about
$(n-1)$-qua\-si\-groups where $n$ is the predicate arity.

By definition, an $n$-qua\-si\-group $q$ in invertible in each place; we will
use the notion $\dot q$ for the inversion in the first place:
$$
\dot q(y, x_2 \ldots , x_n)=z  \eqdf q(z, x_2, \ldots , x_n) = y.
$$

\bnote
1)
The subset of $\A^{n+1}$ corresponding to an $n$-qua\-si\-group predicate is called
a \emph{distance-$2$ MDS code} in the theory of error-correcting codes.
Although such codes themselves cannot
correct errors, they are useful in constructions of codes with larger distance.
2) The $n$-dimensional value array of an $n$-qua\-si\-group is known as a \emph{Latin hypercube}.
\enote

\bdefin[(reducible, irreducible)] An $(n-1)$-qua\-si\-group $M$ is called {\it reducible} ({\it irreducible})
iff it can (cannot)
be represented as
$$M\langle x_1,...,x_n\rangle  \equivvv
K\left\langle q(x_{\eta(1)},...,x_{\eta(j)}),x_{\eta(j+1)},...,x_{\eta(n)}\right\rangle $$
where $K$ and $q$ are $(n-j)$- and $j$-quasigroups,
$\eta:\fromOneTo{n}\to\fromOneTo{n}$ is a permutation, and $2\leq j \leq n-2$.
Note that all binary (as well as $1$-ary and $0$-ary) qua\-si\-groups are irreducible by definition
because $2>n-2$ in this case.
\edefin

\bnote
Defined as above, the reducibility property does not depend on the order of the arguments
of a multary qua\-si\-group.
Often (e.\,g. \cite{Belousov}) by reducibility one means the more strict property,
so-called $(i,j)$-re\-du\-ci\-bi\-li\-ty,
when $\eta=(i,i+1,...,n,1,2,...,i-1)$.
We observe this difference to avoid a misunderstanding.
In our definition, the reducibility
corresponds to the $(i,j,\eta)$-re\-du\-ci\-bi\-li\-ty in {\rm \cite{Glu90}},
where $\eta$ is a permutation.
\enote

\bdefin[(isotopic)]
Two $n$-qua\-si\-groups $Q, Q': \A^n \to \A$ are called {\it isotopic} iff
$$Q\langle x_1,\ldots,x_{n+1}\rangle  \equivvv
Q'\big\langle \rho_1(x_1),\ldots,\rho_{n+1}(x_{n+1})\big\rangle $$
where
$\rho_1,\ldots,\rho_{n+1}:\A \to \A$ are $1$-qua\-si\-groups (i.\,e., permutations).
\edefin
\bdefin[(retract)]
If an $l$-ary predicate $K\langle \cdot\rangle $ is obtained
by fixing $n-l > 0$ arguments in an $(n-1)$-qua\-si\-group predicate $M\langle \cdot\rangle $,
then 
$K$ is, obviously, a well-defined
$(l-1)$-qua\-si\-group; this $(l-1)$-qua\-si\-group is called a \emph{retract} of $M$.
\edefin

Our goal is to prove the following theorem.
\btheorem\label{th:1}
Let $M: \A^{n-1} \to \A$ be an $(n-1)$-qua\-si\-group.
Let $K: \A^{k-1} \to \A$ be a maximal {\rm(}by arity\/{\rm)} irreducible retract of $M$
{\rm(}note that $3\leq k\leq n-1${\rm)}.
Suppose $4\leq k\leq n-3$.
Then
\begin{equation}\label{e:th1}
M\langle \bar z\rangle  \equivvv K\big\langle q^1(\bar z^1),\ldots,q^k(\bar z^k)\big\rangle
\end{equation}
where $\bar z^1,\ldots,\bar z^k$ are nonempty pairwise disjoint collections of variables from
$\bar z$ and $q^1,\ldots,q^k$ are multary quasigroups.
\etheorem

\bcor
If the maximum arity of an irreducible retract of a given $n$-qua\-si\-group belongs
to $\{3,\ldots,n-3\}$, then the $n$-qua\-si\-group is reducible.
\ecor

\bnote Theorem~\ref{th:1} is not much more stronger than its corollary:
indeed,
the decomposition (\ref{e:th1}) exists for every reducible multary quasigroup $M$
and every irreducible retract $K$ that is maximal
in the sense that unfixing one or more
variables always gives a reducible retract.
Such the conclusion can be drawn if we consider a (tree) decomposition
of $M$ into superposition of irreducible multary qua\-si\-groups;
$K$ must be (up to isotopy and changing the order of arguments) an element
of the decomposition.
More results on the structure of decomposition tree of a reducible
multary quasigroup can be found in \cite{Cher}.
\enote

\bnote\label{rem:cases}
1) By numerical reasons \cite{PotKro:asymp}, almost all $n$-qua\-si\-groups
of order $4$ are irreducible with $k=n-1$.

2) If $|\A|\equiv 0\bmod 4$ and $n$ is odd,
then there are irreducible $(n-1)$-qua\-si\-groups
with $k=n-2$ \cite{Kro:n-2}; e.\,g., the $4$-qua\-si\-group with the following value table:
\begin{eqnarray*}
0123\ 1032\ 2310\ 3201
\ \ 1032\ 0123\ 3201\ 2310
\ \ 2301\ 3210\ 1023\ 0132
\ \ 3210\ 2301\ 0132\ 1023
\\[-0.5ex]
1032\ 0123\ 3201\ 2310
\ \ 0123\ 1032\ 2310\ 3201
\ \ 3210\ 2301\ 0132\ 1023
\ \ 2301\ 3210\ 1023\ 0132
\\[-0.5ex]
2310\ 3201\ 0123\ 1032
\ \ 3201\ 2310\ 1032\ 0123
\ \ 0132\ 1023\ 3210\ 2301
\ \ 1023\ 0132\ 2301\ 3210
\\[-0.5ex]
3201\ 2310\ 1032\ 0123
\ \ 2310\ 3201\ 0123\ 1032
\ \ 1023\ 0132\ 2301\ 3210
\ \ 0132\ 1023\ 3210\ 2301
\end{eqnarray*}

3) If $k=3$, or $k=n-2$ and $n$ is odd, or $k=n-2$ and $|\A|\not\equiv 0\bmod 4$,
then the existence of irreducible $n$-qua\-si\-groups is an open question.
\enote

In Section~\ref{s:aux} we consider several simple statements, which will be used later.
Section~\ref{s:prof} is the proof of Theorem~\ref{th:1}, which consists of several steps,
arranged as propositions.  In the Appendix~A we consider the proof of Theorem~\ref{th:1}
by the example of a $6$-qua\-si\-group.
  In the Appendix~B, for convenience, we cite the list of notations.

The author wish to thank the anonymous referees for very helpful suggestions
and for drawing his attention to interesting and useful literature connected with
the subject of this paper.

The results of this paper were announced in \cite{Kro:2004ACCT:decomp}.

\section{Auxiliary statements.}\label{s:aux}

The following two propositions are straightforward.

\blemma\label{p0}
Let $K$ be an $l$-qua\-si\-group and $Q$ be an $(n-l)$-quasigroup.
Then
$$
K\langle \bar x,Q(\bar y)\rangle  \equivvv K(\bar x)=Q(\bar y) \equivvv Q\langle \bar y,K(\bar x)\rangle ,
\qquad \bar x\in \A^l, \  \bar y\in \A^{n-l}.
$$
\elemma

\blemma\label{p00}
Let $M':\A^m\to \A$ be an $m$-qua\-si\-group,
$q$ be a function from $\A^k$ to $\A$, and the predicate $M\langle \cdot\rangle $
is defined by
$$M\langle \bar x, \bar y\rangle \eqdf M'\langle q(\bar x),\bar y\rangle ,\qquad \bar x\in \A^k,\ \bar y\in \A^m.$$
Then $M$ is a well-defined $(k+m-1)$-qua\-si\-group if and only if $q$ is a $k$-quasigroup.
\elemma

The next claim means that a reducible $n$-qua\-si\-group can be represented
as a superposition of retracts.
As a corollary, these retracts uniquely define the multary qua\-si\-group (Lemma~\ref{p000}).

\blemma\label{p0000}
Let $c$ be a $k$-quasigroup, $b$ be an $l$-quasigroup. Let
\begin{equation}\label{eq:f(bc)}
f(\alpha,\bar\beta,\bar\gamma) \defeq c(b(\alpha,\bar\beta),\bar\gamma) ,
\end{equation}
\begin{equation}\label{eq:abc}
c_0(\alpha,\bar\gamma) \defeq f(\alpha,\bar 0,\bar\gamma)
 ,\quad
b_0(\alpha,\bar\beta) \defeq f(\alpha,\bar\beta,\bar 0)
,\quad
a(\alpha) \defeq f(\alpha,\bar 0,\bar 0)
\end{equation}
where $\alpha\in \A$, $\bar\beta\in \A^{l-1}$, $\bar\gamma\in \A^{k-1}$.
Then
\begin{equation}\label{rep}
f(\alpha,\bar\beta,\bar\gamma)\equiv
c_0(a^{-1}(b_0(\alpha,\bar\beta)),\bar\gamma).
\end{equation}
\elemma
\proofr
Substituting (\ref{eq:f(bc)}) to (\ref{eq:abc}) we get
$c_0(\cdot,\bar\gamma)\equiv c(b(\cdot,\bar 0),\bar\gamma)$,
 $b_0(\alpha,\bar\beta)\equiv c(b(\alpha,\bar\beta),\bar 0)$, and
$a(\cdot)\equiv c( b(\cdot,\bar 0),\bar 0)$,
i.\,e.,
$a^{-1}(\cdot)\equiv \dot b(\dot c(\cdot,\bar 0),\bar 0)$.
Using these representations, we can verify the validity of (\ref{rep}):
$$
c_0(a^{-1}(b_0(\alpha,\bar\beta)),\bar\gamma)\equiv
c(b(\dot b(\dot c(c(b(\alpha,\bar\beta),\bar 0),\bar 0),\bar 0),\bar 0),\gamma)\equiv
c(b(\alpha,\bar\beta),\bar\gamma)\equiv
f(\alpha,\bar\beta,\bar\gamma).
$$
\proofend

\blemma\label{p000}
Let $C,$ and $\widetilde C$ be $k$-qua\-si\-groups, $b$ and $\widetilde b$ be
$l$-quasigroups. Suppose
$$ 
C\langle b(\alpha,\bar 0),\bar\gamma,\delta\rangle
\equivvv 
\widetilde C\langle \widetilde b(\alpha,\bar 0),\bar\gamma,\delta\rangle
\ \mbox{ and }\  
C\langle b(\alpha,\bar\beta),\bar 0,\delta\rangle
\equivvv 
\widetilde C\langle \widetilde b(\alpha,\bar\beta),\bar 0,\delta\rangle
$$ 
where $\alpha,\delta\in \A$, $\bar\beta\in \A^{l-1}$, $\bar\gamma\in \A^{k-1}$.
Then $C\langle b(\alpha,\bar\beta),\bar\gamma,\delta\rangle \equivvv \widetilde C\langle \widetilde b(\alpha,\bar\beta),\bar\gamma,\delta\rangle $.
\elemma
%

\section{Theorem proof.}\label{s:prof}

Given $\bar x=(x_1,x_2,...,x_n)$, we use the following notation:
$\bar x^{(k)}\defeq(x_1,...,x_{k-1},x_{k+1},...,x_n)$,
$\bar x^{(k)}\#y\defeq(x_1,...,x_{k-1},y,x_{k+1},...,x_n)$,
and $\bar x^{(l,k)}=\bar x^{(k,l)}\defeq \bar x^{(l)(k)}$ provided $k<l$.

Let $M:\A^{n-1} \to \A$ be an $(n-1)$-qua\-si\-group;
let $K:\A^{k-1} \to \A$ be an irreducible retract of $M$;
and let $k$ be the maximum number for which such retract exists;
for the rest of this section we suppose that $4\leq k\leq n-3$.
Without loss of generality we assume that $K\langle x_1,...,x_k\rangle \equivvv
M\langle x_1,...,x_k,0,...,0\rangle $. Put $m\defeq n-k$, $\bar x \defeq (x_1,...,x_k)$,
$\bar y \defeq (y_1,...,y_m)$.

In the first four propositions we consider the structure of $k$-ary
and $(k-1)$-ary retracts of $M$ with unfixed arguments $x_1,\ldots,x_k$.


\bpro\label{p1}
Let $L_{i;\bar y^{(i)}}\langle \bar x,z\rangle
\eqdf M\langle \bar x,\bar y^{(i)}\#z\rangle $ be a retract of $M$.
Assume that
$K_{\bar y}\langle \bar x\rangle  \eqdf  M\langle \bar x,\bar y\rangle  \equivvv L_{i;\bar y^{(i)}}\langle \bar x,y_i\rangle $
is an irreducible retract of $L_{i;\bar y^{(i)}}$
(here we only suppose but do not yet claim that such a retract exists).
Then $L_{i;\bar y^{(i)}}$ can be represented as
\begin{equation}\label{eq1}
L_{i;\bar y^{(i)}}
\langle x_1,...,x_k,z\rangle
\equivvv
R_{i;\bar y^{(i)}}
\big\langle x_1,...,x_{j-1},q_{i;\bar y^{(i)}}(x_j,z),x_{j+1},...,x_k\big\rangle
\end{equation}
where $j$ depends (essentially or not) on $i$ and ${\bar y^{(i)}}$, i.\,e., $j=j(i,{\bar y^{(i)}})$,
$R_{i;\bar y^{(i)}}$ and $q_{i;\bar y^{(i)}}$
are multary qua\-si\-groups.
\epro
\proofr
The $k$-qua\-si\-group $L_{i;\bar y^{(i)}}$ is reducible because $k<n-1$.
But its retract $K_{\bar y}$ obtained by fixing
the last variable $z:=y_i$ in $L_{i;\bar y^{(i)}}\langle \cdot\rangle $ is irreducible.
So,
in any decomposition of $L_{i;\bar y^{(i)}}\langle \bar x,z\rangle $ the variable $z$
must be grouped with exactly one other variable; i.\,e., $L_{i;\bar y^{(i)}}$ admits
one of the two decompositions
\begin{eqnarray}\label{eq:Rq}
L_{i;\bar y^{(i)}}
\langle x_1,...,x_k,z\rangle
&\equivvv&
R
\big\langle x_1,...,x_{j-1},x_{j+1},...,x_k,q(x_j,z)\big\rangle
\\
\label{eq:qR}
L_{i;\bar y^{(i)}}
\langle x_1,...,x_k,z\rangle
&\equivvv&
Q
\big\langle x_j,z, r(x_1,...,x_{j-1},x_{j+1},...,x_k),\big\rangle
\end{eqnarray}
for some $2$-qua\-si\-group $q$ ($Q$) and $k$-qua\-si\-group $R$ ($r$).
By Lemma \ref{p0}, (\ref{eq:qR}) implies (\ref{eq:Rq}) with $R=r$, $q=Q$.
Permuting the arguments in (\ref{eq:Rq}), we get the representation (\ref{eq1}).
\proofend


\bpro\label{p2}
All the retracts $K_{\bar y}\langle \bar x\rangle  \eqdf M\langle \bar x,\bar y\rangle $, $\bar y\in \A^m$
are pairwise isotopic and thus irreducible$;$ i.\,e.,
\begin{equation}\label{eq:step0}
K_{\bar y}\langle \bar x\rangle  \equivvv K\left\langle \rho^1_{\bar y}(x_1),\ldots,\rho^k_{\bar y}(x_k)\right\rangle
\end{equation}
where $\rho^1_{\bar y}$,\ldots,$\rho^k_{\bar y}$ are permutations $\A\to \A$.
\epro
\proofr We prove the proposition by induction on the number of nonzero elements in $\bar y$.
The base of induction is $K_{\bar 0}\langle \cdot\rangle  \equivvv K\langle \cdot\rangle $.
For the induction step it is sufficient to prove that
\begin{equation} \label{eq2}
K_{\bar y''}\langle \bar x\rangle  \equivvv K_{\bar y'}\big\langle \bar x^{(j)}\#\rho(x_j)\big\rangle
\end{equation}
where $\bar y'=(y_1,...,y_{i-1},0,0,...,0)$, $\bar
y''=(y_1,...,y_{i-1},y_{i},0,...,0)$, $j=j(i,\bar y')\in \fromOneTo{m}$,
$\rho=\rho_{i,\bar y''}$ is a permutation. 
Then, (\ref{eq2}) means that $K_{\bar y'}$ and $K_{\bar y''}$ are isotopic,
and from (\ref{eq:step0}) with $\bar y=\bar y'$
we have (\ref{eq:step0}) with $\bar y=\bar y''$,
where $\rho^j_{\bar y''}=\rho^j_{\bar y'}\rho$ and $\rho^l_{\bar y''}=\rho^l_{\bar y'}$ for all $l\neq j$.

Let us show (\ref{eq2}). Note that $\bar y''^{(i)}=\bar y'^{(i)}=(y_1,...,y_{i-1},0,...,0)$.
By Proposition \ref{p1}
\begin{eqnarray*}
K_{\bar y'}\langle \bar x\rangle
\equivvv&
M\langle \bar x,\bar y'\rangle
&\equivvv
R_{i;\bar y'^{(i)}}
\big\langle x_1,...,x_{j-1},q_{i;\bar y'^{(i)}}(x_{j},0),x_{j+1},...,x_k\big\rangle ,
\\
K_{\bar y''}\langle \bar x\rangle
\equivvv&
M\langle \bar x,\bar y''\rangle
&\equivvv
R_{i;\bar y'^{(i)}}
\big\langle x_1,...,x_{j-1},q_{i;\bar y'^{(i)}}(x_{j},y_i),x_{j+1},...,x_k\big\rangle
\end{eqnarray*}
where $j=j(i,\bar y'^{(i)})$.
We see that (\ref{eq2}) holds with
$\rho(\cdot)=\dot q_{i;\bar y'^{(i)}}(q_{i;\bar y'^{(i)}}(\cdot,y_i),0)$.
\proofend

Our goal is to show that each of the permutations
$\rho^1_{\bar y}$,\ldots,$\rho^k_{\bar y}$
in (\ref{eq:step0})
essentially depends on its own group of parameters from $\bar y$
and these groups are pairwise disjoint.
At the first step (which will be used for an induction step later),
in Propositions \ref{p3} and \ref{p4},
we will prove that for each $i\in[m]$ there exists a representation
like (\ref{eq:step0}) where only one of $\rho^1_{\bar y}$,\ldots,$\rho^k_{\bar y}$
essentially depends on $y_i$. In the final Proposition~\ref{pfin} we will
show (by induction) the existence of such a representation that is common for all $y_i$,
$i\in[m]$.


\bpro\label{p3}
Each $k$-qua\-si\-group
$L_{i;\bar y^{(i)}}\langle \bar x,z\rangle  \eqdf M\langle \bar x,\bar y^{(i)}\#z\rangle $
can be represented
in the form
\begin{equation}\label{eq3}
L_{i;\bar y^{(i)}}\langle x_1,...,x_k,z\rangle
\equivvv
K\big\langle p^1_{i;\bar y^{(i)}}(x_1),\ldots,p^{j-1}_{i;\bar y^{(i)}}(x_{j-1}),
p^{}_{i;\bar y^{(i)}}(x_j,z),
p^{j+1}_{i;\bar y^{(i)}}(x_{j+1}),\ldots,p^k_{i;\bar y^{(i)}}(x_k)\big\rangle
\end{equation}
where $j=j(i,{\bar y^{(i)}})$, \ $p^{}_{i;\bar y^{(i)}}$
is a $2$-quasigroup,
and $p^t_{i;\bar y^{(i)}}$ is a $1$-quasigroup (i.\,e., permutation) for $t\neq j$.
\epro
\proofr
Fixing $z:=0$ in (\ref{eq1}) and applying Proposition \ref{p2},
we find that for each $i$ and $\bar y^{(i)}$ the $(k-1)$-quasigroup $R_{i;\bar y^{(i)}}$
in (\ref{eq1}) is isotopic to $K$.
\proofend


\bpro\label{p4}
In Proposition \ref{p3} the index $j$
does not depend on ${\bar y^{(i)}}$, i.\,e., $j=j(i)$.
\epro
\proofr
Assume the contrary, i.\,e., there exist $i$, $\bar y'^{(i)}$ and $\bar y''^{(i)}$ such that
$j'\defeq j(i,{\bar y'^{(i)}})\neq j''\defeq j(i,{\bar y''^{(i)}})$. Without loss of
generality we can assume that $j'=1$ and $j''=2$.
So,
\begin{eqnarray}
\label{'}
L_{i;\bar y'^{(i)}}\langle x_1,x_2,x_3,...,x_k,z\rangle &\equivvv&
K\big\langle   \phantom{\scriptstyle 1}p^{}
(x_1,z),\,p^{2}
(x_{2})\phantom{,z},\,p^{3}
(x_{3}),\ldots,\,p^{k}
(x_k)\big\rangle , \\
\label{''}
L_{i;\bar y''^{(i)}}\langle x_1,x_2,x_3,...,x_k,z\rangle &\equivvv&
K\big\langle   r^{1}
(x_{1})\phantom{,z},\,\phantom{\scriptstyle 2}r^{}
(x_2,z),\,r^{3}
(x_{3}),\ldots,\,r^{k}
(x_k) \big\rangle .
\end{eqnarray}

The $k$-qua\-si\-group $K'\langle z,x_2,x_3,...,x_k\rangle \eqdf L_{i;\bar
y'^{(i)}}\langle 0,x_2,x_3,...,x_k,z\rangle \equivvv M\langle \bar x^{(1)}\#0,\bar y'^{(i)}\#z\rangle $ is
isotopic to $K$ (see (\ref{'})) and irreducible. By Proposition \ref{p2}
(taking $x_1:=z$) $K'$ is isotopic to
$K''\langle z,x_2,x_3,...,x_k\rangle \eqdf L_{i;\bar
y''^{(i)}}\langle 0,x_2,x_3,...,x_k,z\rangle \equivvv M\langle \bar x^{(1)}\#0,\bar y''^{(i)}\#z\rangle $.
But $K''$ is reducible  because (\ref{''}) gives its decomposition when $x_1=0$
(here we use the condition $k\geq 4$). We get a contradiction. \proofend

Now we see that the function $j(i)$ divides all $y$-variables into $k$ groups,
where each group corresponds to an $x$-variable.
The next proposition is very important; it consider
the structure of a $(k+1)$-ary retract of $M$ with two $y$-variables
that belong to different groups.
This is the only place where we use the condition $k\neq n-2$;
if $k=n-2$, then the proposition does not work,
and $M$ can be irreducible, as noted in Remark~\ref{rem:cases}(2).


\bpro\label{p333}
Let $j(i')=1$, $j(i'')=2$, $v\defeq y_{i'}$, $w\defeq y_{i''}$.
Suppose that values of the variables $\bar y^{(i',i'')}\in \A^{m-2}$
are fixed, and denote by $N(\bar x,v,w)$ the corresponding retract of $M$.
Then
\begin{equation}\label{eq333}
N\langle \bar x,v,w\rangle \equivvv
K\big\langle o^1(x_1,v),
  o^2(x_2,w),
  o^3(x_3),\ldots,
  o^k(x_k)\big\rangle
\end{equation}

%
where
$o^t$, 
$t=1,...,k$ are $2$- and $1$-quasigroups,
which depend on the choice of $i'$, $i''$,  $\bar y^{(i',i'')}$.
\epro
\proofr
Recall that for retracts with variables $v,x_1,x_2,...,x_k$ or $w,x_1,x_2,...,x_k$ we have
the decompositions
\begin{eqnarray}\label{eq:x1v}
&& K\langle p^{}(x_1,v),p^2(x_2),\ldots,p^k(x_k)\rangle,
\\ \label{eq:x2w}
&&  K \langle q^1(x_1),q^{}(x_2,w),q^3(x_3),\ldots,q^k(x_k)\rangle .
\end{eqnarray}
respectively.
Consider possible decompositions of $N$.
Taking into account that fixing $v$ and $w$ results in an irreducible
retract, isotopic to $K$, we can conclude that
$N\langle \bar x,v,w\rangle$
admits one of the following decompositions:
\begin{eqnarray}
 \label{e:adm1}
  N\langle \bar x,v,w\rangle &\equivvv&
  C\langle \bar x,b(v,w) \rangle
  \\
 \label{e:adm2}
  N\langle \bar x,v,w\rangle &\equivvv&
  C\langle \bar x^{(i)}\#b(x_i,v),w \rangle, \qquad i\neq 1
  \\
 \label{e:adm3}
  N\langle \bar x,v,w\rangle &\equivvv&
  C\langle \bar x^{(i)}\#b(x_i,w),v \rangle, \qquad i\neq 2
  \\
 \label{e:adm4}
  N\langle \bar x,v,w\rangle &\equivvv&
  C\langle \bar x^{(i)}\#b(x_i,v,w) \rangle
  \\
 \label{e:adm5}
  N\langle \bar x,v,w\rangle &\equivvv&
  C\langle b(x_1,v), x_2,x_3,\ldots,x_k,w \rangle
  \\
 \label{e:adm6}
  N\langle \bar x,v,w\rangle &\equivvv&
  C\langle x_1,b(x_2,w),x_3,\ldots,x_k,v \rangle
\end{eqnarray}
In the case (\ref{e:adm1}) $C$ must be reducible, and a decomposition of $C$
provides another decomposition of $N$ (in fact, only (\ref{e:adm4}) is suitable).
So, $N$ admits one of (\ref{e:adm2})-(\ref{e:adm6}).
Consider (\ref{e:adm2}).
Fixing $x_1$ and $w$ we get a reducible $k-1$-ary retract with variables $x_2,\ldots,x_k,v$.
But this retract is isotopic to $K$, see (\ref{eq:x1v}), which contradicts to the irreducibility
of $K$. So, (\ref{e:adm2}) is impossible. Similarly, (\ref{e:adm3}) and (\ref{e:adm4}) lead
to contradictions.

Consider (\ref{e:adm5}) (the case (\ref{e:adm6}) is similar).
Again, $C$ must be reducible, and a decomposition of $C$
provides another decomposition of $N$.
Since (\ref{e:adm1})-(\ref{e:adm4}) are inadmissible for $N$,
the only possibility for $C$ is
\[
  C\langle  u,x_2,x_3,\ldots,x_k,w\rangle \equivvv
  C'\langle u,b'(x_2,w),x_3,\ldots,x_k \rangle.
\]
In this case
\[
  N\langle  \bar x,v,w\rangle \equivvv
  C'\langle b(x_1,v),b'(x_2,w),x_3,\ldots,x_k \rangle.
\]
Since $C'$ must be isotopic to $K$, the proposition is proved.
\proofend

Now we are ready to prove the main theorem.
All we need to do is to transform the representation (\ref{eq:step0}) to such a form
that for each $i$ only one of $\rho^1_{\bar y}$, \ldots, $\rho^k_{\bar y}$
(more exactly, only $\rho^{j(i)}_{\bar y}$)
essentially depends on $y_i$.
For induction needs, we formulate a proposition covering all intermediate cases
between Proposition \ref{p1} and Theorem \ref{th:1}.
So, Theorem \ref{th:1} is a partial case of the following proposition,
which will be proved by induction.

Let the function $j:\fromOneTo{m}\to\fromOneTo{k}$ be defined as in Proposition \ref{p4}.
Let
${\mathbf i}^t=\{i^t_1,\ldots,i^t_{m_t}\}$
(where $t\in \fromOneTo{k}$) be the set of all
indexes $i$ such that $j(i)=t$. %
Obviously, $\bigcup_{t=1}^k {\mathbf i}^t=\fromOneTo{m}$ and $\sum_{t=1}^k m_t=m$.
For an arbitrary multiindex ${\mathbf i} =\{i_1,\ldots,i_{m'}\}\subseteq \fromOneTo{m}$ where
$i_1<i_2<\ldots <i_{m'}$ we denote $\bar y_{\mathbf i}\defeq (y_{i_1},\ldots,y_{i_{m'}})$.


\bpro\label{pfin} Let ${\mathbf h}^t\subseteq {\mathbf i}^t$, $t=1,...,k$. Denote ${\mathbf h}=\cup_{t=1}^k {\mathbf h}^t$ and
${\bar {\mathbf h}}=\fromOneTo{m}\backslash {\mathbf h}$. Then for each $\bar y_{\bar {\mathbf h}}$ there
exist $(1+|{\mathbf h}^t|)$-quasigroups $q^t_{\bar y_{\bar {\mathbf h}}}$, $t=1,...,k$ such that
\begin{equation}\label{000}
M\langle \bar x,\bar y\rangle \equivvv
K\big\langle q^1_{\bar y_{\bar {\mathbf h}}}(x_1,\bar y_{{\mathbf h}^1})
,\ldots,q^k_{\bar y_{\bar {\mathbf h}}}(x_k,\bar y_{{\mathbf h}^k})\big\rangle .
\end{equation}
\epro%
\proofr %
Propositions \ref{p3} and \ref{p4} imply that the claim holds for $|{\mathbf h}|=1$.
Let this be the induction base.

Assume the claim holds for $|{\mathbf h}|=b$.
Let us show that it holds for
${\mathbf h}={\mathbf g}\subseteq\fromOneTo{m}$ where $|{\mathbf g}|=b+1$.
We fix arbitrary different $i',i''\in {\mathbf g}$ and denote
${\mathbf d} \defeq {\mathbf g}\backslash \{i',i''\}$,
${\mathbf d}^t = {\mathbf d} \cap {\mathbf i}^t$.
Denote $v\defeq y_{i'}$ and $w\defeq y_{i''}$.
We consider two cases: $j(i')=j(i'')$ and $j(i')\neq j(i'')$.

\underline{Case 1}. Assume $j(i')=j(i'')=1$, without loss of generality.

By the inductive hypothesis for ${\mathbf h}={\mathbf d}\cup\{i'\}$, $\bar {\mathbf h}=\bar {\mathbf g}\cup\{i''\}$, we have

\begin{equation}\label{i'}
M\langle \bar x,\bar y\rangle
\equivvv
K\big\langle p^{}_{w}(x_1,\bar y_{{\mathbf d}^1},v),p^2_{w}(x_2,\bar y_{{\mathbf d}^2}),
\ldots,p^k_{w}(x_k,\bar y_{{\mathbf d}^k})\big\rangle
\end{equation}
where multary quasigroups $p^{}_{w}$, $p^t_{w}$, $t=2,...,k$ depend also on ${\bar y_{\bar {\mathbf g}}}$, %
i.\,e., $p^t_{w}=p^t_{\bar y_{\bar {\mathbf g}},w}$.

By the inductive hypothesis for ${\mathbf h}={\mathbf d}\cup\{i''\}$, $\bar {\mathbf h}=\bar {\mathbf g}\cup\{i'\}$, we have
$$
M\langle \bar x,\bar y\rangle \equivvv
K\big\langle r_{v}(x_1,\bar y_{{\mathbf d}^1},w),r^2_{v}(x_2,\bar y_{{\mathbf d}^2})
\ldots,r^k_{v}(x_k,\bar y_{{\mathbf d}^k})\big\rangle
$$
where multary quasigroups $r^{}_{v}$, $r^t_{v}$, $t=2,...,k$ depend also on ${\bar y_{\bar {\mathbf g}}}$, %
i.\,e., $r^t_{v}=r^t_{\bar y_{\bar {\mathbf g}},v}$.

Equating these two representations of $M$ and setting %
$v:=0$, $\bar y_{{\mathbf d}^1}:=\bar 0$, we obtain
$$
K\big\langle p^{}_{w}(x_1,\bar 0,0),p^2_{w}(x_2,\bar y_{{\mathbf d}^2}),
\ldots,p^k_{w}(x_k,\bar y_{{\mathbf d}^k})\big\rangle
\equivvv
K\big\langle r^{}_{0}(x_1,\bar 0,w),r^2_{0}(x_2,\bar y_{{\mathbf d}^2}),
\ldots,r^k_{0}(x_k,\bar y_{{\mathbf d}^k})\big\rangle .$$ %
Changing the variables as
$u=p^{}_{w}(x_1,\bar 0,0)\,\iff\, x_1=\dot p^{}_{w}(u,\bar 0,0)$,
we get
$$
K\big\langle u,p^2_{w}(x_2,\bar y_{{\mathbf d}^2}),
\ldots,p^k_{w}(x_k,\bar y_{{\mathbf d}^k})\big\rangle
\equivvv
K\big\langle r^{}_{0}(\dot p^{}_{w}(u,\bar 0,0),\bar 0,w),r^2_{0}(x_2,\bar y_{{\mathbf d}^2}),
\ldots,r^k_{0}(x_k,\bar y_{{\mathbf d}^k})\big\rangle .
$$
Substituting $p^{}_{w}(x_1,\bar y_{{\mathbf d}^1},v)$ for $u$, we have
\begin{eqnarray} \nonumber
&&K\big\langle p^{}_{w}(x_1,\bar y_{{\mathbf d}^1},v),p^2_{w}(x_2,\bar y_{{\mathbf d}^2}),
\ldots,p^k_{w}(x_k,\bar y_{{\mathbf d}^k})\big\rangle  \equivvv\\
\nonumber && \equivvv
K\big\langle r^{}_{0}(\dot p^{}_{w}(p^{}_{w}(x_1,\bar y_{{\mathbf d}^1},v),\bar 0,0),\bar 0,w),r^2_{0}(x_2,\bar y_{{\mathbf d}^2}),
\ldots,r^k_{0}(x_k,\bar y_{{\mathbf d}^k})\big\rangle .
\end{eqnarray}
Since, by (\ref{i'}), the left part is equivalent to $M\langle \bar x,\bar y\rangle $, we have (\ref{000})
with ${\mathbf h}={\mathbf g}$, ${\mathbf h}^1={\mathbf d}^1\cup\{i',i''\}$, ${\mathbf h}^t={\mathbf d}^t$ for $t\neq 1$,
$q^1_{\bar y_{\bar {\mathbf h}}}(x_1,\bar y_{{\mathbf h}^1})=r^{}_{0}(\dot p^{}_{w}(p^{}_{w}(x_1,\bar y_{{\mathbf d}^1},v),\bar 0,0),\bar 0,w)$, and
$q^t_{\bar y_{\bar {\mathbf h}}}=r^t_{\bar y_{\bar {\mathbf d}},0}$ for $t\neq 1$. By Lemma
\ref{p00}, the 
function $q^1_{\bar y_{\bar {\mathbf h}}}$ is a
multary qua\-si\-group.

\vspace{0.5em}\underline{Case 2}. Assume $j(i')=1$, $j(i'')=2$, without loss of generality.

By the inductive hypothesis, for every  ${\bar y_{\bar {\mathbf g}}}$ we have
\begin{eqnarray}\nonumber
 M\langle \bar x,\bar y\rangle  &\equivvv&
 K\big\langle
   p^{}_{w}(x_1,\bar y_{{\mathbf d}^1},v),
 \ p^2_{w}(x_2,\bar y_{{\mathbf d}^2}),\hphantom{,w}
 \ p^3_{w}(x_3,\bar y_{{\mathbf d}^3}),\ldots,
 \ p^k_{w}(x_k,\bar y_{{\mathbf d}^k})\big\rangle
\\
M\langle \bar x,\bar y\rangle   &\equivvv&
 K\big\langle
   r^1_{v}(x_1,\bar y_{{\mathbf d}^1}),\hphantom{,v}\,
 \ r^{}_{v}(x_2,\bar y_{{\mathbf d}^2},w),\,
 \ r^3_{v}(x_3,\bar y_{{\mathbf d}^3}),\ldots,\,
 \ r^k_{v}(x_k,\bar y_{{\mathbf d}^k})\big\rangle.
 \label{i''}
\end{eqnarray}
Repeating steps of Case 1, we derive
\begin{equation} \label{eq:eq1case}
M\langle \bar x,\bar y\rangle
\equivvv
K\big\langle s_{w}(x_1,\bar y_{{\mathbf d}^1},v),r^{}_{0}(x_2,\bar y_{{\mathbf d}^2},w),r^3_{0}(x_3,\bar y_{{\mathbf d}^3}),
\ldots,r^k_{0}(x_k,\bar y_{{\mathbf d}^k})\big\rangle
\end{equation}
where $s_{w}(x_1,\bar y_{{\mathbf d}^1},v)\defeq r^1_{0}(\dot p^{}_{w}(p^{}_{w}(x_1,\bar y_{{\mathbf d}^1},v),\bar y_{{\mathbf d}^1},0),\bar y_{{\mathbf d}^1})$.
It remains to eliminate the $w$-dependence of the formula in the first position of $K\langle \ldots\rangle $.
Put
\begin{equation}\label{ghj}
\widetilde M\langle \bar x,\bar y\rangle  \eqdf
K\big\langle s_{0}(x_1,\bar y_{{\mathbf d}^1},v),r^{}_{0}(x_2,\bar y_{{\mathbf d}^2},w),r^3_{0}(x_3,\bar y_{{\mathbf d}^3}),
\ldots,r^k_{0}(x_k,\bar y_{{\mathbf d}^k})\big\rangle .
\end{equation}
Setting $w:=0$ in (\ref{ghj}) and (\ref{eq:eq1case}),
we find that
$\widetilde M\langle \bar x,\bar y^{(i'')}\#0\rangle \equivvv M\langle \bar x,\bar y^{(i'')}\#0\rangle $.
On the other hand, $s_{w}(x_1,\bar y_{{\mathbf d}^1},0)\equiv r^1_{0}(x_1,\bar y_{{\mathbf d}^1})$ by definition of $s_w$;
therefore, setting $v:=0$ in (\ref{ghj}) and (\ref{i''}),
we get
$\widetilde M\langle \bar x,\bar y^{(i')}\#0\rangle \equivvv M\langle \bar x,\bar y^{(i')}\#0\rangle $.
Considering $M$ and $\widetilde M$ as $3$-qua\-si\-groups with the arguments $x_1$, $x_3$, $v$, $w$
and parameters $\bar x^{(1,3)}$, $\bar y^{(i',i'')}$,
and
taking into account the decompositions (\ref{eq333}) and (\ref{ghj}),
we see by Lemma~\ref{p000}
(with $\alpha=x_1 $, $\bar\beta=v $, $\delta=x_3 $, $\bar\gamma=w $)
that $M\langle \bar x,\bar y\rangle \equivvv \widetilde M\langle \bar x,\bar y\rangle $.
\proofend


\section*{Appendix A. An example}

In this appendix we consider the proof of Theorem~\ref{th:1} (Proposition~\ref{pfin})
by the example of a $6$-qua\-si\-group $M$.
Assume that all $5$-ary and $4$-ary retracts of $M$ are reducible;
and assume that the $3$-ary retract
$K\langle \bar x\rangle \eqdf M\langle \bar x,0,0,0\rangle$
is irreducible.
Suppose that some $4$-ary retracts of $M$ admit the following decompositions:
\begin{eqnarray*}
M\langle \bar x,y_1,0,0\rangle&\equivvv& R_1\langle q_1(x_1,y_1),x_2,x_3,x_4\rangle\\
M\langle \bar x,0,y_2,0\rangle&\equivvv& R_2\langle q_2(x_1,y_2),x_2,x_3,x_4\rangle\\
M\langle \bar x,0,0,y_3\rangle&\equivvv& R_3\langle x_1,q_3(x_2,y_3),x_3,x_4\rangle.
\end{eqnarray*}
By Proposition~\ref{p2} 
$$ \forall y_1,y_2,y_3: M\langle \bar x, \bar y\rangle \equivvv
 K\langle \rho^1_{y_1,y_2,y_3}(x_1),\rho^2_{y_1,y_2,y_3}(x_2),\rho^3_{y_1,y_2,y_3}(x_3),\rho^4_{y_1,y_2,y_3}(x_4)\rangle $$
 where $\rho^1_{y_1,y_2,y_3},\rho^2_{y_1,y_2,y_3},\rho^3_{y_1,y_2,y_3},\rho^4_{y_1,y_2,y_3}:\Sigma\to\Sigma$ are permutations ($1$-qua\-si\-groups).
By Propositions~\ref{p3} and~\ref{p4} we also have
\begin{eqnarray}
\label{e:bas1}
\forall y_2,y_3: M\langle \bar x, \bar y\rangle &\equivvv&
 K\langle p^{}_{1;y_2,y_3}(x_1,y_1),p^2_{1;y_2,y_3}(x_2),p^3_{1;y_2,y_3}(x_3),p^4_{1;y_2,y_3}(x_4)\rangle,
\\
\label{e:bas2}
\forall y_1,y_3:  M\langle \bar x, \bar y\rangle &\equivvv&
 K\langle p^{}_{2;y_1,y_3}(x_1,y_2),p^2_{2;y_1,y_3}(x_2),p^3_{2;y_1,y_3}(x_3),p^4_{2;y_1,y_3}(x_4)\rangle,
\\
\label{e:bas3}
\forall y_1,y_2: M\langle \bar x, \bar y\rangle &\equivvv&
 K\langle p^1_{3;y_1,y_2}(x_1),p^{}_{3;y_1,y_2}(x_2,y_3),p^3_{3;y_1,y_2}(x_3),p^4_{3;y_1,y_2}(x_4)\rangle
\end{eqnarray}
 for some $1$-qua\-si\-groups $p^2_{1;y_2,y_3},p^3_{1;y_2,y_3},p^4_{1;y_2,y_3}$,
 $p^2_{2;y_1,y_3},p^3_{2;y_1,y_3},p^4_{2;y_1,y_3}$,
 $p^1_{3;y_1,y_2},p^3_{3;y_1,y_2},p^4_{3;y_1,y_2}$ and $2$-qua\-si\-groups
 $p^{}_{1;y_2,y_3}$, $p^{}_{2;y_1,y_3}$, $p^{}_{3;y_1,y_2}$.
 So, $y_1$, $y_2$ are grouped with $x_1$ and $y_3$ is grouped with $x_2$;
 i.\,e., $j(1)=j(2)=1$, $j(3)=2$, $\mathbf{i}^1=\{1,2\}$, $\mathbf{i}^2=\{3\}$, $\mathbf{i}^3=\emptyset$, $\mathbf{i}^4=\emptyset$.

 By Proposition~\ref{p333} we have
\begin{equation}\label{e:oooa}
\forall x_2, y_2: M\langle \bar x, \bar y\rangle \equivvv
 K\langle o^1_{y_2}(x_1,y_1),o^2_{y_2}(x_2,y_3),o^3_{y_2}(x_3),o^4_{y_2}(x_4)\rangle
\end{equation}
 for some $o^1_{y_2}$, $o^2_{y_2}$, $o^3_{y_2}$, $o^4_{y_2}$.

From (\ref{e:bas1})-(\ref{e:bas3}) we see that Proposition~\ref{pfin}
holds for $\mathbf{h}=\{1\}$, $\mathbf{h}=\{2\}$, and $\mathbf{h}=\{3\}$.

1) We will prove that it holds for $\mathbf{h}=\{1,3\}$.
Let $i'=1$ and $i''=3$. Since $j(i')=1\neq j(i'')=2$,
we have the situation of \underline{Case 2}.
Equating (\ref{e:bas1}) and (\ref{e:bas3}) and setting $y_1:=0$ we obtain
\begin{eqnarray*}
 && \nonumber
 K\big\langle p^{}_{1;y_2,y_3}(x_1,0),p^2_{1;y_2,y_3}(x_2),p^3_{1;y_2,y_3}(x_3),p^4_{1;y_2,y_3}(x_4)\big\rangle
 \\&&\qquad\qquad\qquad\qquad\qquad\qquad \equivvv\label{e:st10}
 K\big\langle p^1_{3;0,y_2}(x_1),p^{}_{3;0,y_2}(x_2,y_3),p^3_{3;0,y_2}(x_3),p^4_{3;0,y_2}(x_4)\big\rangle.
\end{eqnarray*}
Substituting $x_1:=\dot p^{}_{1;y_2,y_3}(u,0)$ we get
\begin{eqnarray*}
 && \nonumber
 K\big\langle u,p^2_{1;y_2,y_3}(x_2),p^3_{1;y_2,y_3}(x_3),p^4_{1;y_2,y_3}(x_4)\big\rangle
 \\&&\qquad\qquad\qquad\qquad\qquad \equivvv\label{e:st11}
 K\big\langle p^1_{3;0,y_2}(\dot p^{}_{1;y_2,y_3}(u,0)),p^{}_{3;0,y_2}(x_2,y_3),p^3_{3;0,y_2}(x_3),p^4_{3;0,y_2}(x_4)\big\rangle.
\end{eqnarray*}
Substituting $u:= p^{}_{1;y_2,y_3}(x_1,y_1)$ we get
\begin{eqnarray}
 && \nonumber
 K\big\langle p^{}_{1;y_2,y_3}(x_1,y_1),p^2_{1;y_2,y_3}(x_2),p^3_{1;y_2,y_3}(x_3),p^4_{1;y_2,y_3}(x_4)\big\rangle
 \\&&\qquad \equivvv\label{e:st12}
 K\big\langle p^1_{3;0,y_2}(\dot p^{}_{1;y_2,y_3}(p^{}_{1;y_2,y_3}(x_1,y_1),0)),p^{}_{3;0,y_2}(x_2,y_3),p^3_{3;0,y_2}(x_3),p^4_{3;0,y_2}(x_4)\big\rangle.
\end{eqnarray}
Since, by (\ref{e:bas1}), the left part of (\ref{e:st12}) is equivalent to $M\langle\bar x,\bar y\rangle$,
we have the following:
\begin{equation*}
\label{eq:eq1case-a}
 M\langle\bar x,\bar y\rangle
 \equivvv
 K\big\langle s_{y_2,y_3}(x_1,y_1),p^{}_{3;0,y_2}(x_2,y_3),p^3_{3;0,y_2}(x_3),p^4_{3;0,y_2}(x_4)\big\rangle
\end{equation*}
where $s_{y_2,y_3}(x_1,y_1)\defeq p^1_{3;0,y_2}(\dot p^{}_{1;y_2,y_3}(p^{}_{1;y_2,y_3}(x_1,y_1),0))$.
To eliminate the subindex $y_3$, define
\begin{equation}
\label{ghj-a}
 \widetilde M\langle\bar x,\bar y\rangle
 \equivvv
 K\big\langle s_{y_2,0}(x_1,y_1),p^{}_{3;0,y_2}(x_2,y_3),p^3_{3;0,y_2}(x_3),p^4_{3;0,y_2}(x_4)\big\rangle.
\end{equation}
It remains to check that $M$ and $\widetilde M$ coincide.
Firstly, $M\langle\bar x,y_1,y_2,0\rangle\equivvv\widetilde M\langle\bar x,y_1,y_2,0\rangle$.
Secondly, from $s_{y_2,y_3}(x_1,0)\equiv p^1_{3;0,y_2}(x_1)$ and (\ref{e:bas3}) we
derive that $M\langle\bar x,0,y_2,y_3\rangle\equivvv\widetilde M\langle\bar x,0,y_2,y_3\rangle$.
For any fixed $x_2$, $x_4$, $y_2$ we have decompositions of both
$M\langle\bar x,\bar y\rangle$ and $\widetilde M\langle\bar x,\bar y\rangle$
of type $C(b(x_1,y_1),y_3,x_3)$,
see (\ref{e:oooa}) and (\ref{ghj-a}). By Lemma~\ref{p000}
$M\langle\bar x,\bar y\rangle\equivvv\widetilde M\langle\bar x,\bar y\rangle$,
and, thus, for some $s^1_{y_2}$, $s^2_{y_2}$, $s^3_{y_2}$, $s^4_{y_2}$ we have
\begin{equation}
\label{bass1}
  M\langle\bar x,\bar y\rangle
 \equivvv
 K\big\langle s^1_{y_2}(x_1,y_1),s^2_{y_2}(x_2,y_3),s^3_{y_2}(x_3),s^4_{y_2}(x_4)\big\rangle.
\end{equation}

2) Similarly, the statement holds for $\mathbf h=\{2,3\}$, and
for some $r^1_{y_1}$, $r^2_{y_1}$, $r^3_{y_1}$, $r^4_{y_1}$ we have
\begin{equation}
\label{bass2}
  M\langle\bar x,\bar y\rangle
 \equivvv
 K\big\langle r^1_{y_1}(x_1,y_2),r^2_{y_1}(x_2,y_3),r^3_{y_1}(x_3),r^4_{y_1}(x_4)\big\rangle.
\end{equation}

3) Now, we are ready to prove the statement for $\mathbf{h}=\{1,2,3\}$.
Let $i'=1$ and $i''=2$. Since $j(i')= j(i'')$,
we have the situation of \underline{Case 1}.
The representations (\ref{bass1}) and (\ref{bass2}) play the role of the induction hypothesis;
equating them and setting $y_1:=0$ we get
\begin{equation*}
 K\big\langle s^1_{y_2}(x_1,0),s^2_{y_2}(x_2,y_3),s^3_{y_2}(x_3),s^4_{y_2}(x_4)\big\rangle
 \equivvv
 K\big\langle r^1_{0}(x_1,y_2),r^2_{0}(x_2,y_3),r^3_{0}(x_3),r^4_{0}(x_4)\big\rangle.
\end{equation*}
Substitute $x_1:=\dot s^1_{y_2}(u,0)$:
\begin{equation*}
 K\big\langle u,s^2_{y_2}(x_2,y_3),s^3_{y_2}(x_3),s^4_{y_2}(x_4)\big\rangle
 \equivvv
 K\big\langle r^1_{0}(\dot s^1_{y_2}(u,0),y_2),r^2_{0}(x_2,y_3),r^3_{0}(x_3),r^4_{0}(x_4)\big\rangle.
\end{equation*}
Substituting $u:=s^1_{y_2}(x_1,y_1)$
and denoting $r(x_1,y_1,y_2)\defeq r^1_{0}(\dot s^1_{y_2}(s^1_{y_2}(x_1,y_1),0),y_2)$,
we obtain
\begin{equation*}
 K\big\langle s^1_{y_2}(x_1,y_1),s^2_{y_2}(x_2,y_3),s^3_{y_2}(x_3),s^4_{y_2}(x_4)\big\rangle
 \equivvv
 K\big\langle r(x_1,y_1,y_2),r^2_{0}(x_2,y_3),r^3_{0}(x_3),r^4_{0}(x_4)\big\rangle.
\end{equation*}
By (\ref{bass1}), the left part is equivalent to $M\langle \bar x, \bar y\rangle$.
Since $r(x_1,y_1,y_2)$ is a $3$-qua\-si\-group, by Lemma~\ref{p00}, Theorem~\ref{th:1}
for our example is proved.

\section*{Appendix B. Notation list}\vspace{-0.5cm}
\begin{itemize}
    \item
      $\A$ is a nonempty set; $\A^n$ is the set of $n$-words over $\A$.
    \item
      $0$ is some fixed element of $\A$; $\bar 0$ is the all-zero word.
    \item
      $\fromOneTo{n} \defeq \{1,\ldots,n\}$.
    \item
      $q\langle x_1, \ldots , x_n, x_{n+1}\rangle  \eqdf q(x_1, \ldots , x_n) = x_{n+1}$.
    \item
      $\dot q(y, x_2 \ldots , x_n)=z  \eqdf q(z, x_2, \ldots , x_n) = y$.
    \item
      If $\bar x=(x_1,x_2,...,x_n)$, then\\
        $\bar x^{(k)}\defeq(x_1,...,x_{k-1},x_{k+1},...,x_n)$, \\
        $\bar x^{(k)}\#y\defeq(x_1,...,x_{k-1},y,x_{k+1},...,x_n)$,\\
        $\bar x^{(l,k)}=\bar x^{(k,l)}\defeq \bar x^{(l)(k)}$ where $k<l$.
\end{itemize}

\vspace{-0.5cm}
\providecommand\href[2]{#2} \providecommand\url[1]{\href{#1}{#1}}

\end{document}